\documentclass[a4paper, 11pt, twoside]{amsart}
\usepackage{srollens-en}[2011/03/03]
\usepackage[left=3.5cm, right=3.5cm, headsep=6mm,
footskip=10mm, top=35mm, bottom=35mm, footnotesep=5mm, headheight=2cm]{geometry}
\usepackage{todonotes}

\usepackage{enumitem}
\newlist{prooflist}{description}{1}
\setlist[prooflist]{font=\normalfont\itshape, labelindent=\parindent, leftmargin=0pt}
\usepackage{soul,xcolor}
\usepackage{mathtools}

\newcommand{\R}{\mathbb{R}}

\newcommand{\C}{\mathbb{C}}

\renewcommand{\bar}{\overline} 
\setstcolor{blue}
\usepackage{comment}

\makeatletter
\g@addto@macro\@floatboxreset\centering
\makeatother

\usepackage[unicode, bookmarks, pdftex]{hyperref}
\hypersetup{colorlinks=true, citecolor=NavyBlue, linkcolor=NavyBlue,
            urlcolor=Orange, pdfpagemode=UseNone, breaklinks=true}

\makeatletter
\def\author@andify{%
  \nxandlist {\unskip ,\penalty-1 \space\ignorespaces}%
    {\unskip {} \@@and~}%
    {\unskip \penalty-2 \space \@@and~}%
}
\makeatother

\title[Biholomorphism type of left-invariant complex structures]{Biholomorphism type of left-invariant complex structures on nilpotent Lie groups}

\author[K.~Hasegawa]{Keizo Hasegawa}
\address{Keizo Hasegawa: Department of Mathematics, Graduate School of Science, Osaka University, Toyonaka, Osaka 560-0043, Japan; and Department of Mathematics, Faculty of Education, Niigata University, Ikarashi-nino-cho, Nishi-ku, Niigata 950-2181, Japan}
\email{hasegawa@math.sci.osaka-u.ac.jp}
\email{hasegawa@ed.niigata-u.ac.jp}

\author[S.~Rollenske]{S\"onke Rollenske}
\address{S\"onke Rollenske: FB~12/Mathematik und Informatik, Philipps-Universit\"at Marburg, Hans-Meerwein-Str.~6, 35032 Marburg, Germany}
\email{rollenske@mathematik.uni-marburg.de}

\author[L.~Sillari]{Lorenzo Sillari}
\address{Lorenzo Sillari: Dipartimento di Scienze Matematiche, Fisiche e Informatiche, Unit\`a di Matematica e Informatica, Universit\`a degli Studi di Parma, Parco Area delle Scienze 53/A, 43124 Parma, Italy}
\email{lorenzo.sillari@unipr.it}

\author[A.~Tomassini]{Adriano Tomassini}
\address{Adriano Tomassini: Dipartimento di Scienze Matematiche, Fisiche e Informatiche, Unit\`a di Matematica e Informatica, Universit\`a degli Studi di Parma, Parco Area delle Scienze 53/A, 43124 Parma, Italy}
\email{adriano.tomassini@unipr.it}

\author[K.~Wehler]{Konstantin Wehler}
\address{Konstantin Wehler: FB~12/Mathematik und Informatik, Philipps-Universit\"at Marburg, Hans-Meerwein-Str.~6, 35032 Marburg, Germany}
\email{konstantin.wehler@uni-marburg.de}

\begin{document}

\setcounter{secnumdepth}{3}
\setcounter{tocdepth}{2}

\begin{abstract}
In this note we prove a conjecture by Hasegawa stating that a simply connected, nilpotent Lie group of dimension $2n$ endowed with a left-invariant complex structure is biholomorphic to $\mathbb{C}^n$.
\end{abstract}

\let\thefootnote\relax\footnotetext{K.H.\ is supported by JSPS (Japan Society for the Promotion of Science) KAKENHI Grant Number 25K07000. S.R.\ and K.W.\ were supported by the DFG through the grant RO~3734/4-1. L.S.\ is partially supported by INdAM - GNSAGA Project (code E53C24001950001). A.T.\ is partially supported by GNSAGA of INdAM. L.S.\ and A.T.\ are partially supported by the Project PRIN 2022 ``Real and Complex Manifolds: Geometry and Holomorphic Dynamics'' (code 2022AP8HZ9).}

\maketitle

\section{Introduction}\label{sec:intro}

The uniformisation theorem is a cornerstone of the theory of complex manifolds. It characterises the topological, algebraic-geometric and differential-geometric properties of a Riemann surface in terms of its universal cover being isomorphic to $\mathbb{P}^1$, $\mathbb{C}$ or $\mathbb{H}$. The problem of characterising manifolds with a fixed type of universal cover or determining properties of the universal cover of a certain type of manifold has continued to see interest ever since, and gained traction especially in K\"ahler and algebraic geometry over the last decades, see for instance \cite{K95, CD14, GKT18}.

In the case of non-K\"ahler manifolds the picture is, as usual, less systematic and focuses on the homogeneous case. Here, the most prominent conjecture has been proposed by Hasegawa \cite{Has09, Has10}.

\begin{conj}\label{conj:biholomorphic}
Let $G$ be a simply connected, nilpotent Lie group of dimension $2n$ and let $J$ be a left-invariant complex structure on $G$. Then the complex manifold $(G,J)$ is biholomorphic to $\mathbb{C}^n$.
\end{conj}

We prove the conjecture by explicitly constructing a biholomorphism, which turns out to be polynomial in exponential coordinates.

\begin{thm}\label{thm:A}
Let $G$ be a simply connected, nilpotent Lie group of dimension $2n$ and let $J$ be a left-invariant complex structure on $G$. Then the complex manifold $(G,J)$ is biholomorphic to $\mathbb{C}^n$ by a polynomial map. In particular, the universal cover of a complex nilmanifold of dimension $n$ is biholomorphic to $\mathbb{C}^n$.
\end{thm}

This theorem was already known in several special cases: Nakamura addressed the case of complex nilpotent Lie groups \cite{Nak75}, while the surface case was implicitly contained in Kodaira's work \cite{K66}, and explicitly in the work of Oeljeklaus--Richthofer \cite{OR84}.

Note that complex nilmanifolds, that is, compact quotients of nilpotent Lie groups with left-invariant complex structure, provide a rich source of examples of compact non-K\"ahler manifolds, and have been intensively studied since Thurston's construction of a compact symplectic manifold admitting no K\"ahler metric, known as the Kodaira--Thurston manifold \cite{Thu76}.

In this context, so-called nilpotent complex structures have been introduced in \cite{CFGU97}, and recently Kanda \cite{Kan26} gave an explicit constructive proof of Theorem~\ref{thm:A} in this special case. In terms of the universal cover, nilpotent complex structures can be characterised by $(G,J)$ being an iterated holomorphic principal bundle for the additive group $\mathbb{C}$, so that this case can also be abstractly deduced from the Oka--Grauert principle \cite{G58}. 

In Example \ref{ex:nonnilpotent} we explicitly compute a polynomial biholomorphism from a $6$-dimensional Lie group with a non-nilpotent complex structure to $\C^3$, exhibiting the generality of our result. 

The key ingredient for the general case, beyond the classical work of Snow \cite{S85, S86}, is a theorem by Bialynicki-Birula and Rosenlicht \cite{BBR62} stating that an injective polynomial self-map of $\R^n$ is already surjective. Surjectivity of injective algebraic self-maps was first established by Newman for polynomials in $\mathbb{R}^2$ \cite{New60}, and by Bialynicki-Birula and Rosenlicht for polynomial self-maps of $\mathbb{R}^n$ and for regular self-maps of complex algebraic varieties \cite{BBR62}. The complex case was also obtained independently by Ax, Grothendieck and Borel \cite{Ax69, Gro66, Bor69}.


As a direct consequence of the biholomorphism constructed in the proof of Theorem \ref{thm:A} being polynomial, the applications obtained by Kanda \cite{Kan26} under the nilpotency assumption on $J$ extend to arbitrary left-invariant complex structures. 
A further application to the moduli of complex nilmanifolds is developed in the thesis of Wehler \cite{W26}.


One should note that the analogue of Theorem~\ref{thm:A} fails if the Lie group $G$ is no longer nilpotent, but merely solvable. This can already be observed for the class of two-dimensional complex solvmanifolds, known as Inoue  surfaces (cf.\ \cite{Has05}), whose universal cover is biholomorphic to $\mathbb{C} \times \mathbb{H}$. Nonetheless, Hasegawa \cite{Has10} conjectured that, for $G$ a simply connected, unimodular, solvable Lie group with a left-invariant complex structure $J$, the complex manifold $(G,J)$ is always Stein.

\subsection*{Acknowledgements}
The proof of Conjecture \ref{conj:biholomorphic} was reached in parallel and independently by K.\,Wehler in his PhD thesis \cite{W26}, supervised by S.\,Rollenske and N.\,Istrati, and by K.\,Hasegawa, L.\,Sillari and A.\,Tomassini. When we learned about the existence of each other's proofs, we decided to publish the result together.

\section{Preliminaries and notation}\label{sec:prelim}

Let $G$ be a simply connected, nilpotent Lie group of dimension $2n$ with associated Lie algebra $\mathfrak{g}$. A left-invariant almost complex structure on $G$ is an endomorphism $J \colon \mathfrak{g} \to \mathfrak{g}$ with $J^2 = -\id$, extended to all of $TG$ by left translations. The almost complex structure $J$ yields a decomposition of the complexified Lie algebra
\[
\mathfrak{g}_{\mathbb{C}} = \mathfrak{g}^{1,0} \oplus \mathfrak{g}^{0,1}
\]
into the $\pm i$-eigenspaces of its complex-linear extension, and $J$ is integrable if and only if $\mathfrak{g}^{1,0}$ is a complex Lie subalgebra of $\mathfrak{g}_{\mathbb{C}}$, that is, if $[\mathfrak{g}^{1,0}, \mathfrak{g}^{1,0}] \subseteq \mathfrak{g}^{1,0}$. From now on we will always assume that $J$ is integrable, in which case $J$ defines a left-invariant complex structure on $G$, and the pair $(G,J)$ is a complex manifold of complex dimension $n$.

In \cite{S86}, Snow constructed a map $\Phi \colon (G,J) \to \mathbb{C}^n$, which is used crucially in the proof of Theorem~\ref{thm:A}. We will briefly recall his construction. Since $G$ is simply connected, we have an inclusion $G \hookrightarrow G_{\mathbb{C}}$ into the universal complexification of $G$ given by the inclusion $\mathfrak{g} \hookrightarrow \mathfrak{g}_{\mathbb{C}}$. Moreover, the exponential map $\exp \colon \mathfrak{g}_{\mathbb{C}} \to G_{\mathbb{C}}$ is a biholomorphism because $G_{\mathbb{C}}$ is also simply connected and nilpotent. In particular, the subgroups $\exp(\mathfrak{g}^{1,0}) = G^{1,0}$ and $\exp(\mathfrak{g}^{0,1}) = G^{0,1}$ of $G_{\mathbb{C}}$ are closed. After choosing a basis of $\mathfrak{g}^{1,0}$, the composition
\begin{equation}\label{eq:quotient is C^n}
  \begin{tikzcd}
    \mathbb{C}^n \cong \mathfrak{g}^{1,0} \arrow[r, "\exp"]
      & G^{1,0} \arrow[hookrightarrow]{r}
      & G_{\mathbb{C}} \arrow[r, "\pi"]
      & G_{\mathbb{C}} / G^{0,1}
  \end{tikzcd}
\end{equation}
defines a biholomorphism $\psi \colon G_{\mathbb{C}}/G^{0,1} \to \mathbb{C}^n$. We denote by $\Phi \colon (G,J) \to \mathbb{C}^n$ the restriction of $\psi \circ \pi \colon G_{\mathbb{C}} \to \mathbb{C}^n$ to $G \subset G_{\mathbb{C}}$, and by \cite[Theorem~1, p.~194]{S86} the map $\Phi$ is a local biholomorphism.

\begin{rem}\label{rem: Snow}
    Snow also proves this result in the case where $G$ is just solvable, see Remark~\ref{rem: bounded domain}. In our setting of nilpotent groups his proof can be somewhat streamlined as follows: we first show that the composition given in \eqref{eq:quotient is C^n} defines a biholomorphism $G_\C/G^{0,1}\cong \C^n$.
Consider the multiplication map
\begin{equation*}\label{eq:mu}
\mu \colon G^{1,0} \times G^{0,1} \to G_{\mathbb{C}}.
\end{equation*}
In exponential coordinates, the decomposition $\mathfrak{g}_{\mathbb{C}} = \mathfrak{g}^{1,0} \oplus \mathfrak{g}^{0,1}$ and the Baker--Campbell--Hausdorff formula show that $\mu$ is a polynomial map $\mathbb{C}^{2n} \to \mathbb{C}^{2n}$. It is injective because $G^{1,0} \cap G^{0,1} = \{e\}$. The Bialynicki-Birula--Rosenlicht Theorem (or the Ax--Grothendieck Theorem) now states that the map $\mu$ is also surjective, and therefore a biholomorphism. Hence, the group $G_{\mathbb{C}}$ is biholomorphic to $G^{1,0} \times G^{0,1}$, and thus $G_{\mathbb{C}}/G^{0,1} \cong \mathbb{C}^n$.

 It remains to show that $\Phi$ is a local biholomorphism. The map $\Phi$ is equivariant with respect to the left action of $G$. Hence, $\Phi$ has constant rank. The differential at the identity $d\Phi_e \colon \mathfrak{g} \to \mathfrak{g}_{\mathbb{C}}/ \mathfrak{g}^{0,1} $ is the natural projection. It is $\C$-linear, so $\Phi$ is holomorphic, and injective since $\ker d\Phi_e = \mathfrak{g} \cap \mathfrak{g}^{0,1} = \{0\}$. As both sides have real dimension $2n$, the map $d\Phi_e$ is an isomorphism, and $\Phi$ is a local biholomorphism.
\end{rem}


\section{Proof of Theorem~\ref{thm:A}}

We will show that the map $\Phi \colon (G,J) \mapsto \IC^n$ is bijective and polynomial in exponential coordinates. It follows directly from the construction of $\Phi$ that it is injective if $G$ is nilpotent. Indeed, two elements in $(G,J)$ have the same image under $\Phi$ if and only if they differ by an element of $G\cap G^{0,1}$, but since the exponential map is bijective we have $G \cap G^{0,1} = \exp(\mathfrak{g} \cap \mathfrak{g}^{0,1}) = \exp(0)$.

In order to show that $\Phi$ is also surjective we have to describe the map more explicitly. Let $Z_1, \dots, Z_n$ be a basis of $\mathfrak{g}^{1,0}$, and identify $\mathbb{C}^{2n}$ with $\mathfrak{g}_{\mathbb{C}}$ via the basis $Z_1, \dots, Z_n$, $\bar{Z}_1, \dots, \bar{Z}_n$. Consider the splitting map
\[
\Psi \colon \mathbb{C}^{2n} \to \mathbb{C}^{2n}, \qquad
\Psi(w, u) = \log \left( \exp \left(\sum_{j=1}^n w_j Z_j\right)
             \exp \left(\sum_{k=1}^n u_k \bar{Z}_k\right) \right),
\]
where $w =(w_1, \dots, w_n), u= (u_1, \dots, u_n) \in \C^n$, and $\log$ denotes the inverse of the exponential map. By the Baker--Campbell--Hausdorff formula, which terminates because $\mathfrak{g}_{\mathbb{C}}$ is nilpotent, $\Psi$ is a polynomial map. Moreover, $\Psi$ is injective because $G^{1,0}\cap G^{0,1} = \{e\}$. Now, the determinant of the Jacobian matrix of $\Psi$ is a complex polynomial $D$ that vanishes nowhere, since $\Psi$ is a local diffeomorphism. By the weak Hilbert's Nullstellensatz \cite[Corollary~7.10, p.~82]{AM69} applied to the ideal $\langle D \rangle$, we have that $D$ is a unit in the ring of complex polynomials, hence a non-zero constant. By \cite[Theorem~2.1, p.~294]{BCW82}, an injective polynomial map $\mathbb{C}^{2n} \to \mathbb{C}^{2n}$ with constant Jacobian is invertible with polynomial inverse. Therefore $\Psi^{-1}$ is polynomial. Finally, every real element $g \in G \subset G_{\mathbb{C}}$ can be written as
\[
    g = \exp \left(\sum_{i=1}^n z_i Z_i + \bar{z}_i \bar{Z}_i\right),
\]
and the factorisation
\[
    g = \exp \left(\sum_{j=1}^n w_j Z_j\right)\exp \left(\sum_{k=1}^n u_k \bar{Z}_k\right)
\]
is given by $(w_1, \dots, w_n, u_1, \dots, u_n) = \Psi^{-1}(z, \bar{z})$. By the construction of $\Phi$ we have $\Phi(g) = (w_1, \dots, w_n)$, so the map $\Phi \circ \exp \colon \mathfrak{g} \to \IC^n$ is polynomial in the variables $z_j, \bar{z}_j$. Identifying $\mathfrak{g}$ with $\mathbb{R}^{2n}$ via the real and imaginary parts of $Z_j$, the map $\Phi \circ \exp \colon \mathbb{R}^{2n} \to \mathbb{C}^n \cong \mathbb{R}^{2n}$ is a real polynomial map. Since we already showed that $\Phi$ is injective, so is $\Phi \circ \exp$. Hence, by the Bialynicki-Birula--Rosenlicht Theorem, $\Phi \circ \exp$ is also surjective. Therefore, $\Phi$ is a biholomorphism, polynomial in exponential coordinates, and $(G, J) \cong \mathbb{C}^n$.


\section{Concluding Remarks}\label{sect: examples}

If the complex structure $J$ is nilpotent in the
sense of Cordero--Fern\'andez--Gray--Ugarte \cite{CFGU97}, then there exists a basis $Z_1, \dots, Z_n$ of $\mathfrak{g}^{1,0}$ compatible with the ascending central series of $\mathfrak{g}$. In this case, it follows directly from the argument in the previous section that the degree of the polynomials given by the map $\Phi$ is bounded by the nilpotency index of the Lie algebra $\mathfrak{g}$. This is explicitly demonstrated by the examples given in \cite[Section~5]{Kan26}.
However, as the following example shows, this bound no longer holds if the complex structure is non-nilpotent.

\begin{exam}\label{ex:nonnilpotent}
The biholomorphism $\Phi$ given in the proof of Theorem~\ref{thm:A} is polynomial for \emph{every} left-invariant complex structure, including those that are non-nilpotent. In real dimension $6$ there are, up to isomorphism, exactly two nilpotent Lie algebras carrying non-nilpotent complex structures \cite{COUV16}; we treat one of them.

Let $\mathfrak{g} = \mathfrak{h}_{19}^-$ be the $3$-step nilpotent Lie algebra with basis $\{X_1, Y_1, X_2, Y_2, X_3, Y_3\}$ and non-trivial brackets
\[
[X_1, X_2] = -Y_3, \quad [X_1, X_3] = -X_2, \quad
[Y_1, Y_2] = -Y_3, \quad [Y_1, X_3] = -Y_2,
\]
endowed with the complex structure $J X_k = Y_k$, for $k = 1, 2, 3$. Its center $\langle Y_3 \rangle$ is one-dimensional and contains no non-zero $J$-invariant subspace, so $J$ is non-nilpotent. Writing $g = \exp\!\big(\sum_{k} x_k X_k + y_k Y_k\big)$ and $z_k = \tfrac{1}{2}(x_k + i y_k)$, a direct computation with the Baker--Campbell--Hausdorff formula gives $\Phi(g) = (w_1, w_2, w_3)$ with
\begin{align*}
  w_1 &= z_1, \\
  w_2 &= z_2 + \tfrac{1}{2} z_1 \bar{z}_3 + \tfrac{i}{4} z_1 \big(\bar{z}_1 z_2 - z_1 \bar{z}_2\big) + \tfrac{i}{12}\, \lvert z_1 \rvert^2 z_1 (z_3 + \bar{z}_3), \\
  w_3 &= z_3 + \tfrac{i}{2}\big(z_1 \bar{z}_2 - \bar{z}_1 z_2\big) - \tfrac{i}{6}\, \lvert z_1 \rvert^2(z_3 + \bar{z}_3).
\end{align*}
This is a polynomial diffeomorphism $\mathbb{R}^6 \to \mathbb{C}^3$ of degree four, one larger than the nilpotency index, illustrating Theorem~\ref{thm:A} in a case beyond the reach of \cite{Kan26}. 
Following the proof of Theorem \ref{thm:A} given in \cite{W26} we obtain that a general bound on the degree of the polynomial map $\Phi$ is given by twice the nilpotency index of $\mathfrak{g}$.
\end{exam}

\begin{rem}\label{rem: bounded domain}
    As already mentioned in Remark \ref{rem: Snow}, the map $\Phi$ can also be defined if $G$ is just solvable. However, its image is in general only biholomorphic to a domain in $\C^n$. We observe that any homogeneous bounded domain in $\C^n$ is a Siegel domain, which has a non-nilpotent solvable Lie group structure.
\end{rem}

\bibliographystyle{alpha}
\bibliography{bibtex}

@book {AM69,
    AUTHOR = {Atiyah, M. F. and Macdonald, I. G.},
     TITLE = {Introduction to commutative algebra},
 PUBLISHER = {Addison-Wesley Publishing Co., Reading, Mass.-London-Don Mills, Ont.},
      YEAR = {1969},
     PAGES = {ix+128},
   MRCLASS = {13.00},
  MRNUMBER = {242802},
MRREVIEWER = {Johnny\ A.\ Johnson},
}

@article {Ax69,
    AUTHOR = {Ax, J.},
     TITLE = {Injective endomorphisms of varieties and schemes},
   JOURNAL = {Pacific J. Math.},
    VOLUME = {31},
      YEAR = {1969},
     PAGES = {1--7},
      ISSN = {0030-8730,1945-5844},
   MRCLASS = {14.15},
  MRNUMBER = {251036},
MRREVIEWER = {M.\ Miyanishi},
       URL = {http://projecteuclid.org/euclid.pjm/1102978046},
}

@article{BBR62,
    AUTHOR = {Bialynicki-Birula, A. and Rosenlicht, M.},
     TITLE = {Injective morphisms of real algebraic varieties},
   JOURNAL = {Proc. Amer. Math. Soc.},
    VOLUME = {13},
      YEAR = {1962},
     PAGES = {200--203},
      ISSN = {0002-9939,1088-6826},
   MRCLASS = {14.15},
  MRNUMBER = {140516},
MRREVIEWER = {P.\ Roquette},
       DOI = {10.2307/2033904},
       URL = {https://doi.org/10.2307/2033904},
}

@article {BCW82,
    AUTHOR = {Bass, H. and Connell, E. H. and Wright, D.},
     TITLE = {The {J}acobian conjecture: reduction of degree and formal expansion of the inverse},
   JOURNAL = {Bull. Amer. Math. Soc. (N.S.)},
    VOLUME = {7},
      YEAR = {1982},
    NUMBER = {2},
     PAGES = {287--330},
      ISSN = {0273-0979,1088-9485},
   MRCLASS = {14H20 (13B25)},
  MRNUMBER = {663785},
MRREVIEWER = {Stephen\ McAdam},
       DOI = {10.1090/S0273-0979-1982-15032-7},
       URL = {https://doi.org/10.1090/S0273-0979-1982-15032-7},
}

@article {Bor69,
    AUTHOR = {Borel, A.},
     TITLE = {Injective endomorphisms of algebraic varieties},
   JOURNAL = {Arch. Math. (Basel)},
    VOLUME = {20},
      YEAR = {1969},
     PAGES = {531--537},
      ISSN = {0003-889X,1420-8938},
   MRCLASS = {14.15},
  MRNUMBER = {255545},
MRREVIEWER = {Andrzej\ Bia\l ynicki-Birula},
       DOI = {10.1007/BF01899460},
       URL = {https://doi.org/10.1007/BF01899460},
}

@article {CD14,
    AUTHOR = {Catanese, F. and Di Scala, A. J.},
     TITLE = {A characterization of varieties whose universal cover is a bounded symmetric domain without ball factors},
   JOURNAL = {Adv. Math.},
    VOLUME = {257},
      YEAR = {2014},
     PAGES = {567--580},
      ISSN = {0001-8708,1090-2082},
   MRCLASS = {32Q30 (14C30 14G35 32J25 32M15 32N05 32Q20)},
  MRNUMBER = {3187659},
MRREVIEWER = {Azniv\ Kasparian},
       DOI = {10.1016/j.aim.2014.02.030},
       URL = {https://doi.org/10.1016/j.aim.2014.02.030},
}

@article{CFGU97,
    AUTHOR = {Cordero, L. A. and Fern\'andez, M. and Gray, A. and Ugarte, L.},
     TITLE = {Nilpotent complex structures on compact nilmanifolds},
 BOOKTITLE = {Proceedings of the Workshop on Differential Geometry and Topology (Palermo, 1996)},
   JOURNAL = {Rend. Circ. Mat. Palermo (2) Suppl.},
    volume = {49},
      YEAR = {1997},
     PAGES = {83--100},
   MRCLASS = {53C15 (53C30 57T15)},
  MRNUMBER = {1491422},
}

@incollection {GKT18,
    AUTHOR = {Greb, D. and Kebekus, S. and Taji, BeB.hrouz},
     TITLE = {Uniformisation of higher-dimensional minimal varieties},
 BOOKTITLE = {Algebraic geometry: {S}alt {L}ake {C}ity 2015},
    SERIES = {Proc. Sympos. Pure Math.},
    VOLUME = {97.1},
     PAGES = {277--308},
 PUBLISHER = {Amer. Math. Soc., Providence, RI},
      YEAR = {2018},
      ISBN = {978-1-4704-3577-6},
   MRCLASS = {32Q30 (14E20 14E30 32Q26)},
  MRNUMBER = {3821153},
MRREVIEWER = {Yu-Chao\ Tu},
       DOI = {10.1090/pspum/097.1/01676},
       URL = {https://doi.org/10.1090/pspum/097.1/01676},
}

@article {G58,
    AUTHOR = {Grauert, H.},
    TITLE = {Analytische {F}aserungen \"uber holomorph-vollst\"andigen {R}\"aumen},
    JOURNAL = {Math. Ann.},
    VOLUME = {135},
    YEAR = {1958},
    PAGES = {263--273},
    ISSN = {0025-5831,1432-1807},
    MRCLASS = {32.00},
    MRNUMBER = {98199},
    MRREVIEWER = {H.\ Cartan},
    DOI = {10.1007/BF01351803},
    URL = {https://doi.org/10.1007/BF01351803},
}

@article{Gro66,
  author  = {Grothendieck, A.},
  title   = {Éléments de géométrie algébrique ({IV}):
             {É}tude locale des schémas et des morphismes
             de schémas ({Troisième partie})},
  journal = {Publ. Math. Inst. Hautes \'Etudes Sci.},
  volume  = {28},
  year    = {1966},
  pages   = {1--255},
}

@article{Has05,
    AUTHOR = {Hasegawa, K.},
     TITLE = {Complex and {K}\"ahler structures on compact solvmanifolds},
    JOURNAL = {J. Symplectic Geom.},
    VOLUME = {3},
      YEAR = {2005},
    NUMBER = {4},
     PAGES = {749--767},
      ISSN = {1527-5256,1540-2347},
   MRCLASS = {32J15 (22E25 32J27)},
  MRNUMBER = {2235860},
MRREVIEWER = {Azniv\ Kasparian},
       DOI = {10.4310/jsg.2005.v3.n4.a9},
       URL = {https://doi.org/10.4310/jsg.2005.v3.n4.a9},
}

@article{Has10,
    AUTHOR = {Hasegawa, K.},
     TITLE = {Small deformations and non-left-invariant complex structures on six-dimensional compact solvmanifolds},
   JOURNAL = {Differential Geom. Appl.},
    VOLUME = {28},
      YEAR = {2010},
    NUMBER = {2},
     PAGES = {220--227},
      ISSN = {0926-2245,1872-6984},
   MRCLASS = {32G05 (32M10 53C30)},
  MRNUMBER = {2594464},
MRREVIEWER = {Azniv\ Kasparian},
       DOI = {10.1016/j.difgeo.2009.10.003},
       URL = {https://doi.org/10.1016/j.difgeo.2009.10.003},
       
}

@incollection{Has09,
    AUTHOR = {Hasegawa, K.},
     TITLE = {Complex and {K}\"ahler structures on compact homogeneous manifolds---their existence, classification and moduli problem},
 BOOKTITLE = {Singularities---{N}iigata--{T}oyama 2007},
    SERIES = {Adv. Stud. Pure Math.},
    VOLUME = {56},
     PAGES = {151--167},
 PUBLISHER = {Math. Soc. Japan},
      YEAR = {2009},
      ISBN = {978-4-931469-55-6},
   MRCLASS = {32M10 (22E25 53C30 53C55)},
  MRNUMBER = {2604081},
MRREVIEWER = {G.\ K.\ Sankaran},
       DOI = {10.2969/aspm/05610151},
       URL = {https://doi.org/10.2969/aspm/05610151},
}

@article{Kan26,
    AUTHOR = {Kanda, S.},
     TITLE = {Holomorphic polynomial crystallographic actions of nilpotent groups},
      YEAR = {2026},
      NOTE = {Preprint, arXiv:2606.04996},
}

@article {K66,
    AUTHOR = {Kodaira, K.},
     TITLE = {On the structure of compact complex analytic surfaces. {II}},
   JOURNAL = {Amer. J. Math.},
    VOLUME = {88},
      YEAR = {1966},
     PAGES = {682--721},
      ISSN = {0002-9327,1080-6377},
   MRCLASS = {32.44 (57.60)},
  MRNUMBER = {205280},
MRREVIEWER = {M.\ F.\ Atiyah},
       DOI = {10.2307/2373150},
       URL = {https://doi.org/10.2307/2373150},
}

@book {K95,
    AUTHOR = {Koll\'ar, J.},
     TITLE = {Shafarevich maps and automorphic forms},
    SERIES = {M. B. Porter Lectures},
 PUBLISHER = {Princeton University Press, Princeton, NJ},
      YEAR = {1995},
     PAGES = {x+201},
      ISBN = {0-691-04381-7},
   MRCLASS = {14E20 (14J10 14J15 32J18 32N10)},
  MRNUMBER = {1341589},
MRREVIEWER = {Kang\ Zuo},
       DOI = {10.1515/9781400864195},
       URL = {https://doi.org/10.1515/9781400864195},
}

@article{Nak75,
    AUTHOR = {Nakamura, I.},
     TITLE = {Complex parallelisable manifolds and their small deformations},
   JOURNAL = {J. Differential Geometry},
    VOLUME = {10},
      YEAR = {1975},
     PAGES = {85--112},
      ISSN = {0022-040X,1945-743X},
   MRCLASS = {32G10 (32M10)},
  MRNUMBER = {393580},
MRREVIEWER = {J.\ L.\ Koszul},
       URL = {http://projecteuclid.org/euclid.jdg/1214432677},
}

@article {New60,
    AUTHOR = {Newman, D. J.},
     TITLE = {One-one polynomial maps},
   JOURNAL = {Proc. Amer. Math. Soc.},
    VOLUME = {11},
      YEAR = {1960},
     PAGES = {867--870},
      ISSN = {0002-9939,1088-6826},
   MRCLASS = {41.15 (12.99)},
  MRNUMBER = {122817},
MRREVIEWER = {M.\ Rosenlicht},
       DOI = {10.2307/2034426},
       URL = {https://doi.org/10.2307/2034426},
}

@article{OR84,
    AUTHOR = {Oeljeklaus, K. and Richthofer, W.},
     TITLE = {Homogeneous complex surfaces},
   JOURNAL = {Math. Ann.},
    VOLUME = {268},
      YEAR = {1984},
    NUMBER = {3},
     PAGES = {273--292},
      ISSN = {0025-5831,1432-1807},
   MRCLASS = {32M10},
  MRNUMBER = {751730},
MRREVIEWER = {B.\ Gilligan},
       DOI = {10.1007/BF01457059},
       URL = {https://doi.org/10.1007/BF01457059},
}

@article {S86,
	AUTHOR = {Snow, D. M.},
	TITLE = {Invariant complex structures on reductive {L}ie groups},
	JOURNAL = {J. Reine Angew. Math.},
	VOLUME = {371},
	YEAR = {1986},
	PAGES = {191--215},
	ISSN = {0075-4102,1435-5345},
	MRCLASS = {32M10 (17B20)},
	MRNUMBER = {859325},
	MRREVIEWER = {Yusuke\ Sakane},
	DOI = {10.1515/crll.1986.371.191},
	URL = {https://doi.org/10.1515/crll.1986.371.191},
}

@article {S85,
	AUTHOR = {Snow, D. M.},
	TITLE = {Stein quotients of connected complex {L}ie groups},
	JOURNAL = {Manuscripta Math.},
	FJOURNAL = {Manuscripta Mathematica},
	VOLUME = {50},
	YEAR = {1985},
	PAGES = {185--214},
	ISSN = {0025-2611,1432-1785},
	MRCLASS = {32M10 (32E10)},
	MRNUMBER = {784143},
	MRREVIEWER = {R.\ R.\ Simha},
	DOI = {10.1007/BF01168831},
	URL = {https://doi.org/10.1007/BF01168831},
}

@article{Thu76,
    AUTHOR = {Thurston, W. P.},
     TITLE = {Some simple examples of symplectic manifolds},
   JOURNAL = {Proc. Amer. Math. Soc.},
    VOLUME = {55},
      YEAR = {1976},
    NUMBER = {2},
     PAGES = {467--468},
      ISSN = {0002-9939,1088-6826},
   MRCLASS = {53C15 (57D15)},
  MRNUMBER = {0402764},
       DOI = {10.2307/2041749},
       URL = {https://doi.org/10.2307/2041749},
}

@misc {W26,
	author =        {Wehler, K.},
	school =        {Philipps University Marburg},
	type =          {PhD Thesis},
	title =         {Moduli spaces for complex nilmanifolds},
	year =          {2026},
    note = {PhD Thesis: https://doi.org/10.17192/openumr/689},
}

@article{COUV16,
    AUTHOR = {Ceballos, Manuel and Otal, Antonio and Ugarte, Luis 
              and Villacampa, Ra\'ul},
     TITLE = {Invariant complex structures on {$6$}-nilmanifolds: 
              classification, {F}r\"olicher spectral sequence and special 
              {H}ermitian metrics},
   JOURNAL = {J. Geom. Anal.},
  FJOURNAL = {Journal of Geometric Analysis},
    VOLUME = {26},
      YEAR = {2016},
    NUMBER = {1},
     PAGES = {252--286},
      ISSN = {1050-6926,1559-002X},
   MRCLASS = {53C15 (17B30 53C30)},
  MRNUMBER = {3441511},
MRREVIEWER = {Anna\ Fino},
       DOI = {10.1007/s12220-014-9548-4},
       URL = {https://doi.org/10.1007/s12220-014-9548-4},
}

\end{document}